\input amstex
\documentstyle{amsppt}
\magnification=\magstep1
\hoffset= 0.25 true in
\vsize= 9 true in
\NoRunningHeads

\define\RR{\bold R^2}
\define\RRR{\bold R^3}

\define\p{^{\prime}}
\define\pp{^{\prime\prime}}
\define\bproper{$\partial$-proper}
\define\bparallel{$\partial$-parallel}
\define\bd{\partial}
\define\halfspace{\bold R^2 \times [0,\infty)}
\define\prodspace{\bold R^2 \times [0,1]}
\define\homeo{homeomorphic}
\define\3m{3-manifold}
\define\tms{3-manifolds}
\define\inc{incompressible}
\define\irr{irreducible}
\define\birr{$\partial$-irreducible}
\define\rirr{$\bold R^2$-irreducible}

\define\ann{anannular}

\define\inte{int \, }
\define\Inte{Int \, }
\define\FF0{{\Cal F}_0}
\define\Fn0{F_{n_0}}
\define\FFn0{\Cal F_{n_0}}
\define\Cn0{C_{n_0}}
\define\Yn0{Y_{n_0}}
\define\PP{\Cal P}
\define\QQ{\Cal Q}
\define\JJ{\Cal J}

\define\al{\alpha}
\define\be{\beta}
\define\ga{\gamma}

\define\la{\lambda}

\define\EE{\Cal E}

\topmatter
\title End Sums of Irreducible Open 3-Manifolds \endtitle
\author Robert Myers \endauthor
\address Department of Mathematics, Oklahoma State University,
Stillwater, OK 74078-1058, USA  \endaddress
\email myersr\@math.okstate.edu \endemail



\endtopmatter
\document


\head 1. Introduction  \endhead

An end sum is a non-compact analogue of a connected sum. Suppose we are 
given two 
connected, oriented $n$-manifolds $M_1$ and $M_2$. Recall that to form 
their connected sum one chooses an $n$-ball in each $M_i$, removes its 
interior, and then glues together the two $S^{n-1}$ boundary components 
thus created by an orientation reversing homeomorphism. Now suppose 
that $M_1$ and $M_2$ are also open, i.e. non-compact with empty boundary. 
To form an end sum of $M_1$ and $M_2$ one chooses a halfspace $H_i$ 
(a manifold \homeo\ to ${\bold R}^{n-1} \times [0, \infty)$) embedded in 
$M_i$, removes its interior, and then glues together the two resulting 
${\bold R}^{n-1}$ boundary components by an orientation 
reversing homeomorphism. In order for this space $M$ to be an 
$n$-manifold one requires that each $H_i$ be {\bf end-proper} in $M_i$ in 
the sense that its intersection with each compact subset of $M_i$ is 
compact. Note that one can regard $H_i$ as a regular neighborhood of 
an end-proper ray (a 1-manifold \homeo\ to $[0,\infty)$) $\ga_i$ in $M_i$. 

The concept of an end sum was introduced by Gompf in \cite{1}, where he developed 
a smooth version of it to use in the study of exotic $\bold R^4$'s. 
He showed that it induces a monoid structure on the set of oriented 
diffeomorphism classes of smooth 4-manifolds which are \homeo\ to 
$\bold R^4$. He proved that it is well defined by showing that all 
end-proper rays in such a 4-manifold are ambient isotopic. In general, 
the choice of $\ga_i$ determines an end of $M_i$. One can informally regard 
the process of forming an end sum as gluing together an end of 
$M_1$ and an end of $M_2$. Simple examples show that different 
choices of these ends may yield non-\homeo\ $n$-manifolds. However, Gompf's  
arguments generalize to show that for $n \geq 4$ the result depends only on the choice 
of the end-proper homotopy classes of the rays $\ga_i$. (An end-proper map is 
a map under which preimages of compact sets are compact; an end-proper 
homotopy is a homotopy which is an end-proper map.)  

This paper examines the case $n=3$. It compares the resulting theory with 
both that of connected sums of 3-manifolds and of higher dimensional end sums.  
Recall that in general 
the connected sum of 3-manifolds depends on the orientations of the 
summands. The standard examples are the pairs of connected sums of certain 
lens spaces with themselves obtained by choosing the opposite orientation 
on one of the summands. (See \cite{3} or \cite{4}.) A similar phenomenon occurs 
for end sums. Given any two connected, oriented, one-ended, \irr\ open 3-manifolds 
one can choose 
the halfspaces so that a change of orientation in one of the summands 
yields a non-\homeo\ \3m (Theorem 5.1(b)). One can also ensure that neither 
end sum admits an orientation reversing homeomorphism. 
In particular this can be done when both manifolds are \homeo\ to $\RRR$. 

Unlike the higher dimensional theory it turns out that end sums of 
3-manifolds 
do not depend merely on the choice of end-proper homotopy classes of rays. 
In fact, given any two open 3-manifolds as above one can form their 
end sum in uncountably many non-\homeo\ ways along rays in any fixed 
end-proper homotopy classes (Theorem 5.1(c)). 
Again note  
that this applies in particular to end sums of $\RRR$ with itself.

The dependence of an end sum on the choice of halfspaces raises 
difficulties for any attempt to inductively define an end sum with more 
than two summands. We take a different approach to defining such multiple 
end sums which has the advantage of allowing infinitely many summands. 
It is patterned after Scott's work \cite{7} on infinite connected sums of 
3-manifolds. First suppose that we have a countable tree $\Gamma$ to 
each vertex $v_i$ of which we have associated a connected, oriented, 
non-compact \3m $V_i$ whose boundary is a non-empty disjoint union of 
planes. Suppose further that whenever $v_i$ and $v_j$ are joined by an 
edge we have associated to that edge a boundary plane of $V_i$ and a 
boundary plane of $V_j$. The result of gluing each pair of associated 
planes together via an orientation reversing homeomorphism is a 
connected, oriented, non-compact \3m $M$ such that $\bd M$ is either 
empty or a disjoint union of planes. $M$ is called the {\bf plane sum} 
of the $V_i$ along $\Gamma$. 

Now suppose that we are given a countable tree $\Gamma$ to each vertex 
$v_i$ of which we have associated a connected, oriented, open \3m $M_i$. 
Suppose further that to each edge of $\Gamma$ incident with $v_i$ we 
have associated an end of $M_i$ and have chosen an end-proper halfspace 
determining that end. The same end may be associated to different edges, 
but we assume that the halfspaces associated to each edge are distinct 
and disjoint and that their union is end-proper in $M_i$. We then remove 
the interiors of these halfspaces from each $M_i$ to obtain a \3m $V_i$ 
and form the plane sum of these $V_i$ as above. The \3m $M$ so obtained 
is called an {\bf end sum} of the $M_i$ along $\Gamma$. 

In the theory of connected sums of \tms\ primeness and irreducibility 
play a key role, and so one would like appropriate analogues for the 
theory of end sums. Recall that $M$ is prime if whenever it is a 
connected sum of two \tms\ one of the summands must be $S^3$; it is 
\irr\ if every 2-sphere in $M$ bounds a 3-ball in $M$. Irreducible 
\tms\ are prime; every prime, oriented \3m is either \irr\ or \homeo\  
to $S^1 \times S^2$. 

Given the dependence of an end sum on the choice of halfspace and the 
fact that one can obtain non-trivial end sums all of whose summands are 
$\RRR$ it seems reasonable to define an open \3m $M$ to be {\bf end-prime} 
if whenever it is an end sum of two \tms\ one of the manifolds must be 
$\RRR$ and the halfspace in this $\RRR$ must be the standard halfspace 
$\RR \times [0,\infty)$. We say that a non-compact \3m\ $M$ is  $\RR$-{\bf \irr} 
if it is \irr\ and 
every end-proper plane in $M$ bounds an end-proper halfspace in $M$. 
Then it is clear that open $\RR$-\irr\ \tms\ are end-prime. However there 
are uncountably many connected, orientable, \irr, end-prime \tms\  
which fail to be \rirr\ (Theorem 5.8). 

Every compact, oriented \3m\ is either prime or admits a decomposition, 
unique up to homeomorphism, as a connected sum of prime \tms. (See \cite{3} 
or \cite{4}.) 
Scott gave an example \cite{7} of an open \3m which is not prime and is not a 
connected sum, finite or infinite, of prime \tms. His example is 
simply connected but has uncountably many ends. In this paper we give 
an example, inspired by that of Scott, of an open \3m which is not 
end-prime and is not an end sum, 
finite or infinite, of end-prime \tms. Moreover, it is eventually 
end-\irr\ and, unlike Scott's example, is \irr\ and contractible 
(hence is one-ended). See Theorem 6.2. 

The family of 2-spheres arising in a connected sum $M$ has certain 
non-triviality properties. One of the summing 2-spheres bounds a 
3-ball in $M$ if and only if some component of the \3m $M\p$ obtained 
by splitting $M$ along all the summing 2-spheres is a 3-ball. Two 
summing 2-spheres are parallel in $M$ if and only if some  
component of $M\p$ either is $S^2 \times [0, 1]$  or becomes 
$S^2 \times [0, 1]$ upon adding some 3-balls in $M$ bounded by some of its 
boundary 2-spheres. Thus it is very easy to characterize such ``degenerate'' 
connected sums. For end sums and plane sums we regard the sum as being 
{\bf degenerate} if a summing plane is {\bf trivial} (bounds an 
end-proper $\RR \times [0,\infty)$ in $M$), or two such planes are {\bf parallel}  
(cobound an end-proper $\prodspace$ in $M$), or, in the case of plane 
sums, a summing plane is {\bf \bparallel\ } (cobounds an end-proper 
$\prodspace$ with a plane in $\bd M$.) 

Degenerate finite plane sums of \irr\ \tms\ are easily characterized. 
They must have an $\halfspace$ or $\prodspace$ summand (Theorem 
3.1.) For infinite plane sums the situation is more delicate, and 
degenerate sums can arise in surprising ways. However, these can also 
be completely characterized (Theorem 3.2.) Moreover the 
degeneracies can sometimes be eliminated. Any \irr\ open \3m which is 
an end sum of end-prime \tms\ is either end-prime or can be expressed 
as a non-degenerate end sum of end-prime \tms\ (Theorem 6.1.) This is 
a key step in proving that certain \tms\ do not have end-prime 
decompositions. 

The bad behavior of some of our examples when considered as end sums is 
detectable because of their good behavior when considered as plane sums. 
In particular the sums in these examples are non-degenerate, and so the 
summing planes 
are non-trivial and non-parallel. Moreover, they have the very strong 
property that every proper non-trivial plane is ambient isotopic to a 
summing plane. We give some general criteria on the plane 
summands (``strong aplanarity'' and ``anannularity at infinity'') 
which ensure that we obtain such a ``strong'' sum (Theorem 4.1). 
Define a \3m\ $V$ to be {\bf aplanar} if every proper plane in $V$ is 
either trivial or \bparallel. This result 
enables us to prove that if a \3m\ $M$ has a decomposition as a strong 
plane sum, then every decomposition of $M$ as a 
non-degenerate plane sum of aplanar \tms\ along a locally finite tree 
is ambient isotopic to the given strong plane sum (Theorem 4.3). This 
has implications for the mapping class group of $M$ (Corollary 4.4).  

\head 2. Preliminaries \endhead

In this section we state some basic defininitions and some technical 
results from \cite{6}. We then investigate these properties for 
the cases of $\halfspace$ and $\prodspace$. 

We shall work throughout in the PL category. An $m$-manifold $M$ may or
may not have boundary but is assumed to be second countable.
$\bd M$ and $\inte M$ denote the manifold theoretic boundary and interior
of $M$, respectively. Let $A$ be a subset of $M$. The topological boundary,
interior, and closure of $A$ in $M$ are denoted by $Fr_M A$, $Int_M A$, and
$Cl_M A$, respectively, with the subscript deleted when $M$ is clear
from the context. All isotopies of $A$ in $M$ will be ambient.
$A$ is {\bf bounded} if $Cl \, A$ is compact. $M$ is
{\bf open} if $\bd M=\emptyset$ and no component of $M$ is compact.

A map $f:M \rightarrow N$ of manifolds is {\bf \bproper\ } if
$f^{-1}(\bd N)=\bd M$. It is {\bf end-proper} if preimages of compact
sets are compact. It is {\bf proper} if it has both these properties.
These terms are applied to a submanifold if its inclusion map has the
corresponding property.

An {\bf exhaustion} $C=\{C_n\}$ for a  connected, non-compact $m$-manifold
$M$ is a sequence
$C_0 \subseteq C_1 \subseteq C_2 \subseteq \cdots$ of compact, connected 
$m$-submanifolds of $M$ whose union is $M$ such that 
$C_n\cap\bd M$ is either empty or an $(m-1)$-manifold,
$C_n \subseteq Int \, C_{n+1}$, and $M-C_n$ has no bounded
components. Connected non-compact $m$-manifolds always have exhaustions.
A sequence
$V_0 \supseteq V_1 \supseteq V_2 \supseteq \cdots$ of open subsets of
$M$ is an {\bf end sequence} associated to $C$ if each $V_n$ is
a component of $M-C_n$. Two end sequences $\{V_n\}$ and $\{W_p\}$
associated to exhaustions $C$ and $K$ for $M$
are {\bf cofinal} if for every $n$ there is a $p$ such that
$V_n \supseteq W_p$ and for every $p$ there is a $q$ such that
$W_p \supseteq V_q$. Cofinality is an equivalence relation on
end sequences of $M$. The equivalence classes are called the
{\bf ends} of $M$. The set of all ends of $M$ is denoted by
$\varepsilon (M)$. An end-proper map $M \rightarrow N$ induces a
well defined function $\varepsilon(M) \rightarrow \varepsilon(N)$.
If $\bd M$ has no compact components, then the inclusion map induces
a well defined bijection $\varepsilon(int \, M) \rightarrow \varepsilon(M)$.

We refer to \cite{3} and \cite{4} for basic \3m\ topology, including the 
definitions of irreducibility, incompressibility, parallel, \bparallel, 
etc., and of splitting a manifold along a codimension one submanifold. 

We shall need two technical results from \cite{6}. 

\proclaim{Proposition 2.1} Let $M$ be a connected, \irr, non-compact
\3m\  which is not \homeo\  to $\RRR$. Let $\Cal P$ and $\Cal Q$ be
proper surfaces in $M$ which are in general position. Let $\JJ$ be a union
of simple closed curve components of $\PP \cap \QQ$. Assume that the
following conditions are satisfied.
\roster
\item No component of $\Cal P$ or of $\Cal Q$ is a 2-sphere.
\item Each component $J$ of $\Cal J$ bounds
a disk $D(J)$ on $\Cal P$ and a disk $G(J)$ on $\Cal Q$.
\item There is no infinite sequence $\{J_m\}$ of distinct components of
$\Cal J$ such that either $D(J_m) \subseteq int \, D(J_{m+1})$ for
all $m$ or $G(J_m) \subseteq int \, G(J_{m+1})$ for all $m$, i.e.
there is no {\bf infinite nesting} on $\Cal P$ or on $\Cal Q$ among
the components of $\Cal J$.
\endroster
Then there is an
ambient isotopy of $\Cal P$ in $M$, fixed on $\bd M$, which takes $\Cal P$
to a surface $\Cal P\p$ such that $\Cal P\p$ and $\Cal Q$ are in general
position and $(\Cal P\p \cap \Cal Q) \subseteq (\Cal P \cap \Cal Q)-\Cal J$.
Moreover, the isotopy is fixed on $\Cal P\p \cap \Cal Q$.
\endproclaim

\demo{Proof} This is Proposition 2.1 of \cite{6}. \qed \enddemo

\proclaim{Lemma 2.2} Let $M$ be a connected, \irr, non-compact \3m\.
A prop\-er plane $P$ in $M$ is trivial if and only if there
exist sequences $\{D_n\}$ and $\{D_n\p\}$ of disks in $M$ such that $\{D_n\}$
is an exhaustion for $P$, $D_n\p \cap P=\bd D_n$, and $\cup D_n\p$ is
end-proper in $M$.
\endproclaim

\demo{Proof} This is Lemma 1.1 (1) of \cite{6}. \qed \enddemo

\proclaim{Lemma 2.3} $\RRR$ is \rirr. \endproclaim 

\demo{Proof} This follows from the result of Harrold and Moise \cite{2} that 
a topologically 
embedded 2-sphere in $S^3$ which is wild at at most one point bounds a 
3-ball on at least one side. \qed \enddemo

\proclaim{Lemma 2.4} $\halfspace$ is \rirr. \endproclaim

\demo{Proof} Let $\{B_n\}$ be an exhaustion for $\RR$ be concentric disks. 
Let $C_n=B_n \times [0,n+1]$ and $F_n=\bd C_n - \inte B_n$. Suppose $P$ 
is a proper plane in $\halfspace$. Let $\PP=P$ and $\QQ=\cup F_n$. 
Since $P$ is proper and $\{C_n\}$ is an exhaustion there must be an 
infinite sequence $\{K_i\}$ of components of $\PP \cap \QQ$ which is 
nested on $P$. We may assume this sequence is maximal in the sense that if 
$J$ is a component of $\PP \cap \QQ$ which is not in the sequence then 
$D(J)$ does not contain any of the $K_i$. Let $\JJ=(\PP \cap \QQ)-\cup 
K_i$. Apply Proposition 2.1 to eliminate $\JJ$ from the intersection. 
$\PP \cap \QQ$ is now an infinite subsequence of $\{K_i\}$. By 
passing to a further subsequence we may assume that the disks $D\p_i$ 
on $\QQ$ bounded by the $K_i$ are disjoint. Lemma 2.2 now implies that 
$P$ is trivial. \qed \enddemo

It should be noted that $\halfspace$ contains proper planes which are 
not \bparallel. See the discussion in section 3. 

\proclaim{Lemma 2.5} $\prodspace$ is aplanar. \endproclaim

\demo{Proof} Let $\{B_n\}$ be as in Lemma 2.4. We identify $B_n$ with 
$B_n \times \{0\}$. Let $C_n=B_n \times [0,1]$ and $F_n=(\bd B_n) 
\times [0,1]$. 

Suppose $P$ can be isotoped off $C_0$. Then for each $n>0$ we may assume 
that either $P \cap F_n$ is empty or consists of simple closed curves 
which bound disks on $F_n$. It then follows from Lemma 2.2 that $P$ is 
trivial. 

So assume this cannot be done. Let $D$ be a disk in $P$ containing $P \cap 
C_0$ in its interior. Choose $n_0>0$ such that $D \subseteq \Inte C_{n_0}$. 
Let $J_0$ be the component of $P \cap F_{n_0}$ such that $D \subseteq 
\Inte D(J_0)$ and $J_0$ is innermost on $P$ among All such curves. Suppose 
$J$ is a component of $D(J_0) \cap F_{n_0}$ other than $J_0$ which is 
innermost on $P$ among such curves. Then $D(J) \cap F_{n_0}=J$, 
$D(J) \cap D = \emptyset$, and so $D(J)$ lies in $(\RR-\inte B_0) \times 
[0,1]$. Since $F_{n_0}$ is incompressible in this \irr\ \3m\ there is 
an isotopy of $D(J_0)$, fixed on $J_0$, which removes $J$ from 
$D(J_0) \cap F_{n_0}$ and adds no new components. Continuing in this 
fashion we get $D(J_0) \cap F_{n_0}=J_0$. In a similar way we can isotop 
$P$ so as to remove all those components $J$ of $P \cap F_{n_0}$ such that 
$D(J)$ does not contain $D(J_0)$. 

Now let $M=(\RR-\inte B_{n_0}) \times [0,1]$, $\PP=P \cap M$, and 
$\QQ= \cup_{n>n_0} F_n$. Then $\PP$ consists of a half-cylinder 
$S^1 \times [0,\infty)$ and possibly some annuli. These surfaces are 
incompressible in $M$ since otherwise the incompressibility of $F_{n_0}$ 
would imply that $D(J_0)$ could be isotoped off $C_0$. We apply Proposition 
2.1 to isotop $\PP$ so that afterwards $\PP \cap \QQ$ consists of simple 
closed cureves which do not bound disks on $\PP$ or $\QQ$. It follows that 
these curves are concentric on $P$ about $J_0$. Denote them by $J_k$, 
$k \geq 1$, where $D(J_k) \subseteq \inte D(J_{k+1})$. 

Choose $k_0 > 0$ such that $(P \cap C_{n_0}) \subseteq D(J_{k_0})$ and 
$n_1 > n_0$ such that $D(J_{k_0}) \subseteq \inte C_{n_1}$. Next choose 
$k_1 > k_0$ such that $(P \cap C_{n_1}) \subseteq \inte D(J_{k_1})$ and 
$n_2 > n_1$ such that $D(J_{k_1}) \subseteq \inte C_{n_2}$. 
Continuing in this fashion we define sequences $\{k_i\}$ and $\{n_i\}$ 
such that $(P \cap C_{n_i}) \subseteq \inte D(J_{k_i})$ and 
$D(J_{k_i}) \subseteq \inte C_{n_{i+1}}$. We assume that $k_i$ and $n_i$ 
are chosen to be the minimal such indices satisfying these conditions. 
Then $J_{k_i} \subseteq F_{n_i}$. 

For each $i\geq 0$ we have that $P \cap C_{n_i}$ consists of a disk 
$D_i$ with $D(J_{k_{i-1}}) \subseteq D_i \subseteq D(J_{k_i})$ and 
possibly a finite number of annuli. For $i=0$ these annuli lie in 
$C_{n_0}-C_0$ and are each parallel in $C_{n_0}-C_0$ to an annulus 
in $F_{n_0}$ which misses $D_0=D(J_0)$. There is an isotopy of 
$P \cap (C_{n_1} - \Inte C_0)$ in $C_{n_1} - \Inte C_0$, fixed on 
$\bd (C_{n_1} - \Inte C_0)$, which carries it to a surface whose 
intersection with $F_{n_0}$ is $\bd D_0$. Similarly for each even 
$k>0$ these annuli lie in $C_{n_i}-C_{n_{i-1}}$, each of them is 
parallel in $C_{n_i}-C_{n_{i-1}}$ to an annulus in $F_{n_i}$ which 
misses $D_i$, and there is an isotopy of 
$P \cap (C_{n_i}-\Inte C_{n_{i-1}})$ in $C_{n_i}-\Inte C_{n_{i-1}}$, 
fixed on $\bd(C_{n_i}-\Inte C_{n_{i-1}})$, which carries it to a surface 
whose intersection with $F_{n_i}$ is $\bd D_i$. Since these isotopies 
have disjoint compact supports they define an ambient isotopy of $P$ in 
$\prodspace$ after which $P \cap F_{2p}$ 
is a single simple closed curve $K_{2p}$ for each $p \geq 0$, and these 
curves are nested on $P$. 

We may assume $D(K_0)=D_0$. For $p>0$ let $A_{2p}=(B_{n_{2p}}-\inte 
B_{n_{2p-2}})$ and $A\p_{2p}=D(K_{2p})-\inte D(D_{2p-2})$. 
For $p \geq 0$ let $G_{2p}$ be the annulus in $F_{n_{2p}}$ joining 
$\bd B_{2p}$ and $K_{2p}$. Since for $p>0$ we have $C_{n_{2p}}-
\Inte C_{n_{2p-2}}=A_{2p} \times [0,1]$, $A\p_{2p}$ is incompressible 
in $A_{2p} \times [0,1]$, and $A\p_{2p} \cap (A_{2p} \times \{1\})=
\emptyset$ it follows that $A\p_{2p}$ is parallel in $A_{2p} \times 
[0,1]$ to $G_{2p} \cup A_{2p} \cup G_{2p-2}$. It follows that there 
is an embedding of $A_{2p} \times [0,1]$ in $\prodspace$ with 
$A_{2p} \times \{0\}=A_{2p}$, $(\bd A_{2p}) \times [0,1]=G_{2p} 
\cup G_{2p-2}$, and $A_{2p} \times \{1\}=A\p_{2p}$. There is also an 
embedding of $B_0 \times [0,1]$ in $\prodspace$ with $B_{n_0} \times 
\{0\}=B_{n_0}$, $(\bd B_{n_0}) \times [0,1]=G_0$, and $B_{n_0} \times 
\{1\}=D_0$. These embeddings can be chosen so as to agree on the 
$G_{2p}$ and so define a parallelism from $P\p$ to ${\bold R}^2 
\times \{0\}$. \qed \enddemo 

\head 3. Degenerate Plane Sums \endhead

A plane sum $M$ of \tms\ $V_i$ along a tree $\Gamma$ is 
{\bf degenerate} if either (1) some summing plane $E_j$ 
is trivial in $M$, (2) some summing plane $E_j$ is \bparallel\  
in $M$, or (3) some pair of distinct summing planes $E_j$ and 
$E_k$ are parallel in $M$. In this section we give some 
conditions on the $V_i$ which ensure that the plane sum is 
non-degenerate. 

There are some obvious ways to obtain degenerate plane sums, 
such as having some summands \homeo\ to $\halfspace$ or 
$\prodspace$. We shall see below that for a finite plane sum to 
be degenerate such summands must be present. For infinite plane 
sums there are some sources of degeneracy which are only slightly 
less obvious. For example, one can use summands which are \homeo\ 
to a closed 3-ball minus the complement of a disjoint union of open 
disks in its boundary to build a degenerate plane sum having no 
$\halfspace$ or $\prodspace$ summands. 

Other sources of degeneracy are less obvious. We briefly describe 
one such example. A non-compact \3m\ is {\bf almost compact} if it 
is \homeo\ to a compact \3m\ minus a closed subset of its boundary. 
In Example 1 of \cite{8} Scott and Tucker construct 
a \3m\ with interior \homeo\ to $\RRR$ and boundary \homeo\ to 
${\bold R}^2 \times \{0,1\}$ which is not almost compact, hence is 
not \homeo\ to $\prodspace$, 
even though the complement of either boundary plane is \homeo\ to 
$\RR \times [0,1)$. This example has an exhaustion by cylinders $C_n=D_n 
\times [0,1]$, where $C_n$ is embedded in $C_{n+1}$ as a regular 
neighborhood of a knotted arc joining the center of 
$D_{n+1}\times \{0\}$ to that of $D_{n+1}\times \{1\}$. 
Call this \3m\ $V$ and its boundary planes $E$ and 
$E\p$. The plane sum of $V$ and $\halfspace$ obtained by identifying 
$E\p$ with ${\bold R}^2 \times \{0\}$ is \homeo\ to $V-E\p$ and 
hence is \homeo\ to $\halfspace$. If one now takes copies 
$(V_i, E_i, E\p_i)$ of $(V, E, E\p)$ for $i\geq 0$ and forms a 
plane sum $W$ by identifying $E\p_i$ with $E_{i+1}$, then given 
any compact subset $K$ of $W$, there is an embedding of $\halfspace$ 
in $W$ which takes ${\bold R}^2 \times \{0\}$ to $E$ and whose image 
contains $K$. It follows that $W$ is \homeo\ to $\halfspace$. One 
can then take a plane sum of $W$ with an appropriate \3m\ to obtain 
an infinite plane sum $M$ which has a trivial summing plane even though 
there are no $\halfspace$ or $\prodspace$ summands. In particular, 
extending the construction of $W$ to $i<0$ expresses $\RRR$ as such 
a plane sum. (The author thanks B. Winters for first pointing out 
this example to him.) Examples similar to $V$ can be constructed 
having more than two boundary components (\cite{9} or \cite{10}). These can be 
used to build plane sums which have no almost compact summands but 
have summing planes which violate conditions (2) or (3). 

\proclaim{Theorem 3.1} A finite plane sum $M$ of \irr\ \tms\ is degenerate 
if and only if it has an $\halfspace$ or $\prodspace$ summand. 
\endproclaim

\demo{Proof} (1) Suppose $E_j$ is trivial in $M$. Then some component 
of $M-E_j$ has closure $H$ which is homeomorphic to $\prodspace$ 
with $E_j=\RR \times \{0\}$. There is a summand $V_i$ contained in 
$H$ such that $\partial V_i$ is a single plane. By Lemma 2.4  
$\partial V_i$ is trivial in $H$, and it follows that $V_i$ is 
homeomorphic to $\halfspace$. 

(2) Suppose some $E_j$ is \bparallel\ in $M$. Then $M-E_j$ has a 
component whose closure $Q$ is homeomorphic to $\prodspace$, with 
$E_j=\RR \times \{0\}$ and $\RR \times \{1\}$ a component of 
$\partial M$. There must be a summand $V_i$ contained in $Q$ 
which either has exactly one boundary component or has exactly 
two boundary components, one of which is a component of 
$\partial M$. In the first case it follows from Lemma 2.4 that $V_i$ 
is homeomorphic to $\halfspace$; in the second case it follows 
from Lemma 2.5 that $V_i$ is homeomorphic to $\prodspace$. 

(3) Suppose $E_j$ and $E_k$, $j \neq k$, are parallel in $M$. Then 
$M-(E_j \cup E_k)$ has three components, one of which has closure $Q$ 
which is homeomorphic to $\prodspace$, with $E_j=\RR \times \{0\}$ 
and $E_k=\RR \times \{1\}$. Let $M^{\prime}$ be the union of $Q$ with 
the other component whose closure contains $E_j$. Then $E_j$ is 
\bparallel\ in $M^{\prime}$ and the desired conclusion follows 
from (2).  \qed \enddemo 

Suppose $M$ is a plane sum along a tree $\Gamma$ of \irr\ \tms\ $V_i$. 
It will be convenient to adjoin to $\Gamma$ an edge for each 
component of $\bd M$; such edges have one vertex associated to 
a summand and another vertex which is not associated to a summand.  
We call such edges {\bf boundary edges} and the regular edges 
{\bf interior edges}. If a vertex $V_i$ meets an edge $E_j$ such that 
$E_j \cup \inte V_i$ is \homeo\ to $\halfspace$, then we say that 
$V_i$ is a {\bf bad summand} with {\bf bad boundary plane} $E_j$. 
If there are distinct edges $E_j$ and $E_k$ meeting $V_i$ such that 
$E_j \cup E_k \cup \inte V_i$ is \homeo\ to $\prodspace$, then 
$V_i$ is {\bf doubly bad} with a {\bf bad pair} $(E_j, E_k)$ of 
boundary planes. 

A {\bf branch} $\be$ of $\Gamma$ is one of the two components of the 
graph obtained by deleting an interior edge $E_j$ of $\Gamma$. 
We associate to each vertex $V_k$ of $\be$ a {\bf leading edge} 
$E_{\ell}$ as follows. If $V_k$ meets $E_j$, then $E_{\ell}=E_j$. 
If $V_k$ does not meet $E_j$, then $E_{\ell}$ is the unique edge 
meeting $V_k$ which separates $V_k$ from $E_j$. A branch  
$\be$ is {\bf bad} if every vertex $V_i$ of $\be$ is a bad 
summand whose leading edge is a bad boundary plane of $V_i$. 

A {\bf trail} $\al$ of $\Gamma$ joining vertices $V_i$ and $V_j$ 
is the unique reduced edge path between them. If one chooses edges 
$E_p$ and $E_q$ meeting $V_i$ and $V_j$, respectively, which are 
not edges of the trail and one also chooses an orientation for $\al$ 
such that $V_i$ and $V_j$ are, respectively, the first and last 
vertex of $\al$, then to each vertex $V_k$ of $\al$ we associate 
a {\bf leading edge} $E_{\ell}$ and a {\bf lagging edge} $E_m$ as 
follows. If $k=i$, then $E_{\ell}=E_p$ and $E_m$ is the first edge 
in the trail. If $k=j$, then $E_{\ell}$ is the last edge in the trail  
and $E_m=E_q$. If $k \neq i$, $j$, then the trail enters $V_k$ through 
$E_{\ell}$ and exits through $E_m$. The other edges of $\Gamma$ 
meeting $V_k$ are called {\bf side edges}. Each of them determines 
a unique {\bf side branch} which does not contain $V_k$. A trail  
$\al$ is {\bf bad} if every vertex of $\al$ is a doubly bad summand 
whose leading and lagging edges are a bad pair of boundary planes 
and every side branch of $\al$ is bad. 

\proclaim{Theorem 3.2} Let $M$ be a plane sum of \irr\ \tms\ along 
a tree $\Gamma$. 
\roster
\item A summing plane $E_j$ is trivial in $M$ if and only if one 
of the branches determined by deleting $E_j$ from $\Gamma$ is a 
bad branch which contains no boundary edges. 
\item A summing plane $E_j$ is parallel to a boundary component $E$ 
of $M$ if and only if $E_j$ is the leading edge of a vertex $V_p$ 
which is joined by a bad trail to a vertex $V_q$ having $E$ as its 
lagging edge and none of the side branches contains a boundary edge. 
\item Two distinct summing planes $E_j$ and $E_k$ are parallel in 
$M$ if and only if $E_j$ is the leading edge of a vertex $V_p$ which 
is joined by a bad trail to a vertex $V_q$ having $E_k$ as its lagging 
edge and none of the side branches contains a boundary edge. 
\endroster \endproclaim 

\proclaim{Corollary 3.3} A plane sum of \irr\ \tms\ having no 
bad summands is non-degenerate. \qed \endproclaim

\demo{Proof} In each case the sufficiency of the conditions is clear, 
so we prove only their necessity.   

(1) Suppose $E_j$ is trivial in $M$. Then $M-E_j$ has 
a component whose closure $G$ is \homeo\ to $\halfspace$ with 
$E_j={\bold R}^2 \times \{0\}$. The corresponding branch clearly 
has no boundary edges. Let $E_{\ell}$ be the leading edge of a 
vertex $V_k$ of the corresponding branch. By Lemma 2.2 $E_{\ell}$ 
is trivial in $G$ and so $G-E_{\ell}$ has a component whose 
closure $H$ is \homeo\ to $\halfspace$ with 
$E_{\ell}={\bold R}^2 \times \{0\}$. Let $X=E_{\ell} \cup \inte V_k$. 
Suppose $K$ is a compact, connected subset of $X$ such that 
$K \cap E_{\ell} \neq \emptyset$. Then there is a 3-ball $B$ in 
$H$ such that $B \cap E_{\ell}$ is a disk and $K$ lies in $Int_HB$. 
Let $D$ be the closure of $(\bd B)-(B \cap E_{\ell})$. Isotop $D$ in 
$H-K$ so that it is in general position with respect to the union of 
all the summing planes other than $E_{\ell}$. Let $D\p$ be an 
innermost disk on $D$ bounded by one of the components of the 
intersection. Then $D\p$ lies in some $V_r$ and $\bd D\p = \bd D\pp$ 
for a disk $D\pp$ in a component $E_s$ of $\bd V_r$. By the 
irreducibility of $V_r$ we have that $D\p \cup D\pp$ bounds a 3-ball 
$B\p$ in $V_r$. Since $K$ is connected, lies in $V_k$, and meets 
$E_{\ell}$ one has that $B\p \cap K = \emptyset$. Thus an isotopy 
fixed on $K$ can be performed to reduce the number of intersection 
components. Continuing in this fashion there is an isotopy fixed 
on $K$ which carries $B$ into $V_k$. It follows that $X$ is \homeo\ 
to $\halfspace$ and thus $V_k$ is a bad summand with bad boundary 
plane $E_{\ell}$. 

(2) Suppose $E_j$ is parallel in $M$ to a component $E$ of $\bd M$. 
Then $M-E_j$ has a component whose closure $Y$ is \homeo\ to 
$\halfspace$ with $E_j={\bold R}^2 \times \{0\}$ and 
$E={\bold R}^2 \times \{1\}$. Let $V_p$ and $V_q$ be the vertices of 
the corresponding branch determined by $E_j$ such that $E_j$ is in 
$\bd V_p$ and $E$ is in $\bd V_q$. There is a unique trail $\al$ 
joining $V_p$ and $V_q$. 

Suppose $\be$ is a side branch of $\al$ determined by the edge $E_t$. 
Since $E_t$ does not separate $E_j$ from $E$ it follows from Lemma 2.4 
that $E_t$ is trivial in $Y$. Hence by part (1) $\be$ is a bad branch 
which contains no boundary edges. Now suppose $V_k$ is a vertex of 
$\al$ with leading edge $E_{\ell}$ and lagging edge $E_m$. These 
two planes cannot be trivial in $Y$, and so by Lemma 2.4 they must 
be \bparallel\ in $Y$. It follows that they must be parallel to each 
other. Since any other components of $\bd V_k$ determine side edges 
they must be trivial in $Y$, from which it follows that $V_k$ is a 
doubly bad summand with bad pair $(E_{\ell}, E_m)$ of boundary planes. 
Thus $\al$ is a bad trail. 

(3) Suppose $E_j$ and $E_k$, $j\neq k$, are parallel in $M$. Then 
we apply part (2) to the obvious manifold $M\p \subseteq M$ 
such that $E_k$ is a component of $\bd M\p$. \qed \enddemo

\head 4. Strong Plane Sums \endhead

Suppose $V$ is a connected, non-compact, \irr\ \3m\ whose boundary 
is either empty or has each component a plane. 
A {\bf partial plane} is a non-compact, simply connected 2-manifold 
with non-empty boundary. 
$V$ is {\bf strongly aplanar} if it is aplanar and has the property 
that given any proper surface $\PP$ in $V$ each component of which 
is a partial plane, there exists a collar on $\bd V$ containing $\PP$. 
$V$ is {\bf \ann\ at infinity} if for every compact subset $K$ of 
$V$ there is a compact subset $L$ of $V$ containing $K$ such that 
$V-L$ is \ann, i.e. every proper incompressible annulus in $V-L$ 
is \bparallel. We emphasize that in this definition one takes 
$V-L$, not the closure of the complement of a regular neighborhood 
of $L$. 

By Lemma 2.4 we have that $\halfspace$ is \rirr\ and hence aplanar; it is clearly 
\ann\ at infinity. It is not strongly aplanar: Let $\PP=\cup P_n$, where 
$P_n=\{(x,y,z) \, | \, x^2+z^2=n^2\}$. By Lemma 2.5 we have that 
$\prodspace$ is aplanar. It is not strongly aplanar: Let $\PP=
\{(x,0,z) \, | \, -\infty < x < \infty, \, 0 \leq z \leq 1\}$. It is 
also not \ann\ at infinity since any compact subset $L$ is contained 
in a 3-ball of the form $D \times [0,1]$ for some disk $D$ in $\RR$, 
and $(\bd D) \times [0,1]$ is not \bparallel\ in the complement of $L$. 

A {\bf strong plane sum} is a non-degenerate plane sum of \irr, 
strongly aplanar \tms\ each of which is \ann\ at infinity. In this 
section we prove that a strong plane sum has the property that all 
of its expressions as a non-degenerate plane sum of aplanar \tms\ 
along a locally finite tree are unique up to ambient isotopy. We 
treat a slightly more general situation (which arises in the next 
section) by allowing non-separating planes. 

\proclaim{Theorem 4.1} Let $M$ be a connected, \irr, non-compact 
\3m\ whose boundary is either empty or has each component a plane. 
Let $\EE$ be a proper surface in $M$ each component of which is 
a plane. Suppose no component of $\EE$ is trivial or 
\bparallel\ and that no two distinct components are parallel in $M$. 
Suppose each component $V_i$ of the manifold $M\p$ obtained by 
splitting $M$ along $\EE$ is strongly aplanar and \ann\ at infinity. 
Then any proper plane $P$ in $M$ which is neither trivial nor 
\bparallel\ in $M$ is ambient isotopic to a component of $\EE$ via 
an isotopy fixed on $\bd M$. \endproclaim

\demo{Proof} Put $P$ in general position with respect to $\EE$. 
If $P \cap \EE = \emptyset$, then $P$ lies in some $V_i$ and we 
are done, so assume the intersection is non-empty. Let $P\p$ 
be the surface obtained by splitting $P$ along $P \cap \EE$. We 
shall denote the two planes in $\bd M\p$ which are identified to 
obtain $E_j$ by $E\p_j$ and $E\pp_j$. 

{\it Case 1}: There is no infinite nesting on $P$ among 
the components of $P \cap \EE$. 

Suppose that there is infinite nesting on $\EE$ among the components 
of $P \cap \EE$. Then there is infinite nesting on some component $E_j$ 
of $\EE$ among the components of $P \cap E_j$. Let $\{\al_n\}$ be a 
maximal nested sequence on $E_j$ of components of $P \cap E_j$. Since 
there is no infinite nesting on $P$ we may pass to a subsequence whose 
elements bound disjoint disks on $P$. Since $P$ is proper it then 
follows from Lemma 2.2 that $E_j$ is trivial in $M$, a contradiction. 
Thus this situation cannot occur. 

By Proposition 2.1 we can remove all the compact components of $P \cap 
\EE$. Then each component $P\p_k$ of $P\p$ is a partial plane. Suppose 
$P$ meets the component $E_j$ of $\EE$. Then $E_j$ lies in $V_i \cap 
V_m$ for some components $V_i$ and $V_m$ of $M\p$, where possibly $i=m$. 
We may assume that $E\p_j \subseteq V_i$ and $E\pp_j \subseteq V_m$. 
Since $V_i$ is strongly aplanar the union of the $P_k\p$ it contains 
must lie in a collar on $\bd V_i$; a similar statement holds for those 
$P_k\p$ contained in $V_m$. Thus a $P\p_k$ cannot meet distinct 
components of $\bd V_i$ or distinct components of $\bd V_m$. It follows 
that $E_j$ is the only component of $\EE$ meeting $P$. Thus 
$P$ lies in $V_i \cup V_m$ and in fact must lie in a regular 
neighborhood of $E_j$. Thus $P$ can be isotoped off $\EE$, and 
we are done. 

{\it Case 2}: There is infinite nesting on $P$ among the 
components of $P \cap \EE$.

Choose a maximal nested sequence $\{\al_n\}$ on $P$ of components of 
$P \cap \EE$. Let $\JJ$ be the union of the remaining components. 
If there is infinite nesting on some component $E_j$ of $\EE$ among 
the components of $\JJ$, then as in the previous case we may pass to 
a subsequence and  
apply Lemma 2.2 to conclude that $E_j$ is trivial in $M$, 
a contradiction. We can therefore apply Proposition 2.1 to eliminate 
$\JJ$ from $P \cap \EE$. If this now has only finitely many components 
we may use irreducibility to perform a finite sequence of isotopies 
which pushes $P$ off $\EE$, and we are done. So assume that $P \cap \EE$ 
is now a nested infinite sequence $\{\al_n\}$. Let $\al\p_n$, $\al\pp_n$ 
denote the preimages of $\al_n$ in $E\p_j$, $E\pp_j$ respectively. 

\proclaim{Lemma 4.2} There is an $N \geq 0$ and a $j$ such that 
for all $n \geq N$ one has that $\al_n$ lies in $E_j$. Moreover, 
$\{\al_n\}_{n \geq N}$ can be re-indexed so as to form a nested sequence 
on $E_j$. If $E\p_j$ and $E\pp_j$ lie in the same component $V_i$ of 
$M\p$, then no component of $P\p$ with boundary in $\cup_{n \geq N} 
(\al\p_n \cup \al\pp_n)$ meets both $E\p_j$ and $E\pp_j$. \endproclaim 

\demo{Proof} If $P$ meets infinitely many $E_j$ then choosing an 
innermost $\al_n$ on each of these $E_j$ yields an end-proper disjoint 
union of disks to which we may apply Lemma 2.2 to conclude that $P$ 
is trivial. Therefore $P$ meets only finitely many components of $\EE$. 

Consider a plane $E_j$ which meets $P$ in infinitely many components. 
Suppose infinitely many of these components bound disks on $E_j$ whose 
interiors miss $P$. Then, as above, there is a subsequence of $\{\al_n\}$ 
which one may use to contradict the non-triviality of $P$. Thus there 
are only finitely many such components, and so after deleting finitely 
many curves the rest can be renumbered so as to form a nested sequence 
on $E_j$. 

Suppose $P$ meets each of two distinct planes $E_j$ and $E_k$ infinitely 
often. We may assume that $E\p_j$ and $E\p_k$ are both boundary 
components of some $V_i$. Then all but finitely many components of 
$P\p \cap V_i$ are annuli. From these we may obtain a sequence $\{A_m\}$ 
of annuli each of which joins $E\p_j$ to $E\p_k$. The union of $A_m$ 
with the disks in $E\p_j$ and $E\p_k$  bounded by the components of 
$\bd A_m$ is a 2-sphere which bounds a 3-ball $B_m$ in $V_i$. Since 
$P$ is proper $\{B_m\}$ is, after renumbering, an exhaustion of $V_i$. 
Let $K=B_0$. Then for every compact subset $L$ of $V_i$ containing 
$K$ there is a $t > 0$ such that $L \subseteq \Inte B_t$. Then 
$A_t$ is a proper incompressible annulus in $V_i-L$. Since $A_t$ joins 
two different components of $\bd(V_i-L)$ it cannot be \bparallel\ in 
$V_i-L$. This contradicts the assumption that $V_i$ is \ann\ at 
infinity, and thus $P$ meets only one component of $\EE$, say $E_j$, 
infinitely often. 

A similar argument proves the assertion about $E\p_j$ and $E\pp_j$ 
lying in the same $V_i$. \qed \enddemo 

Let $V_i$ be the component of $M\p$ containing $E\p_j$. Now 
$P - \inte D_N$ when split along its intersection with $\EE$ 
meets $V_i$ in a family of proper annuli whose boundaries form a 
nested sequence on $E\p_j$ (and on $E\pp_j$ if it also lies in 
$\bd V_i$; in this case none of these annuli meet both $E\p_j$ 
and $E\pp_j$.) 

We construct a sequence $\{A_m\}$ of certain of these annuli with 
$\bd A_m$ in $E\p_j$ as follows. Let $\beta_0$ 
be the innermost of the $\alpha_n$, $n\geq N$. Let $A_0$ be the annulus 
having $\beta_0$ as one boundary component; let $\gamma_0$ be the other 
boundary component. Then $\beta_0 \cup \gamma_0 = \bd A_0\p$ for an annulus 
$A_0\p$ in $E_j$. Suppose $A_0, \cdots ,A_m$ and 
$A_0\p, \cdots ,A_m\p$ have been defined. Let $\beta_{m+1}$ be the 
innermost of the $\alpha_n$, $n \geq N$, which is not contained in 
$A_0\p \cup \cdots \cup A_m\p$. Let $A_{m+1}$ be the annulus having 
$\beta_{m+1}$ as one boundary component; let $\gamma_{m+1}$ be the other 
boundary component. Then $\beta_{m+1} \cup \gamma_{m+1}=\bd A_{m+1}\p$ 
for an annulus $A_{m+1}\p$ in $E_j$. 

Consider the torus $A_m \cup A_m\p$. Compression of this torus along 
the disk in $E_j$ bounded by $\beta_m$ yields a 2-sphere which, by the 
irreducibility of $V_i$, bounds a 3-ball $B_m$ in $V_i$. The non-compactness 
and propriety of $E_j$ imply that $B_m$ contains the compressing disk. 
It follows that $A_m \cup A_m\p$ bounds a compact 3-manifold $Q_m$ 
which either is a solid torus across which $A_m$ and $ A_m\p$ are 
parallel or is \homeo\  to the exterior of a non-trivial knot in $S^3$ for 
which $\beta_m$ is a meridian curve. By construction the $Q_m$ are 
disjoint.

Suppose there are infinitely many $Q_m$ which are \homeo\  to non-trivial 
knot exteriors. Let $K$ be a compact, connected subset of $V_i$ which meets 
the interior of the disk bounded by $\beta_0$. Since $V_i$ is anannular 
at infinity there is a compact subset $L$ of $V_i$ such that $K \subseteq L$ 
and every \bproper\  incompressible annulus in $V_i-L$ is \bparallel. 
Since $P$ is proper in $W$ there are only finitely many $Q_m$ which meet $L$. 
So there is a $Q_p$ which misses $L$ and is \homeo\  to a non-trivial knot 
exterior. The fact that $Q_p$ is not a solid torus implies that $A_p$ is 
not \bparallel\  and so must be compressible in $V_i-L$. Let $D$ be a 
compressing disk, and let $D\p$ be the disk on $E_j$ bounded by $\beta_p$. 
Then there is an annulus $A$ in $A_p$ such that $\bd A=\bd D \cup \bd D\p$. 
The 2-sphere $D \cup D\p \cup A$ bounds a 3-ball $B$ such that 
$B \cap \bd V_i =D\p$. Moreover, $K \subseteq L \subseteq Int(B)$. This 
shows that $V_i$ is \homeo\  to $\halfspace$. As this contradicts the 
non-triviality of $E_j$ it follows that this situation cannot occur. 

We may now assume, by choosing a larger $N$ and re-indexing, that all of the 
$Q_m$ are solid tori across which $A_m$ is parallel to $A_m\p$. Perform an 
ambient isotopy with support in the union of the $Q_m$ to remove all 
$\alpha_n$, $n > N$, from $P \cap E_j$. Then $P \cap \EE$ has only 
finitely many components. They can all be removed by irreducibility in 
the standard way, putting us in the case $P \cap \EE=\emptyset$ of 
Case 1, and we are done. \qed \enddemo 

\proclaim{Theorem 4.3} Let $N$ be a non-degenerate plane sum 
of \irr, aplanar \tms\ along a locally finite tree. Let $M$ 
be a strong plane sum. Let $\PP$ and $\EE$ be the unions of the 
respective sets of summing planes. Suppose $g:N \rightarrow M$ is 
a homeomorphism. Then $g$ is ambient isotopic rel $\bd N$ to a 
homeomorphism $h$ such that $h(\PP)=\EE$. \endproclaim

\demo{Proof} Choose a component $P_0$ of $\PP$. Use Theorem 4.1 to 
isotop $g$ so that $g(P_0)$ is a component, say $E_0$, of $\EE$. 
Now $N-P_0$ has two components. Denote their closures by $X_0$ and 
$Y_0$. Let $W_0$ be the summand of $N$ such that $W_0 \subseteq X_0$ 
and $P_0 \subseteq \bd W_0$. Use Theorem 4.1 finitely many times to 
perform an ambient isotopy of $g$ fixed on $Y_0$ after which 
$g(\PP \cap W_0) \subseteq \EE$. It then follows from the aplanarity 
of $W_0$ and the non-degeneracy of $M$ that $g(W_0)$ is a summand, 
say $V_0$, of $M$. Let $X_1$ be the closure of $X_0-V_0$ in $X_0$. 
Apply the same argument to each component of $X_1$. Continue in this fashion with the successively 
smaller manifolds $X_{n+1}$ obtained by deleting one summand from each 
component of $X_n$ as above. Note that for a given summand there are at most 
finitely many isotopies which are not fixed on that summand, and so one can 
use the sequence of isotopies to define a single ambient isotopy defined on  
$X_0$. Then repeat this procedure on $Y_0$ to finally obtain the isotopy to 
the desired homeomorphism $h$. \qed \enddemo 

\proclaim{Corollary 4.4} Let $M$ be a strong plane sum along a locally 
finite tree; let $\EE$ be the union of the set of summing planes. Suppose 
$g:M \rightarrow M$ is a homeomorphism. Then $g$ is isotopic rel $\bd M$ 
to a homeomorphism $h$ such that $h(\EE)=\EE$. \qed \endproclaim 

\head 5. Strong End Sums \endhead

An end sum $M$ of \tms\ $M_i$ is {\bf strong} if the halfspaces 
in the $M_i$ have been chosen so that the corresponding plane sum is 
strong. In this section we give conditions under which this can be done. 
We also use similar techniques to construct uncountably many end-prime 
\tms\ which are not \rirr. Recall that the choice of an end-proper 
halfspace in $M_i$ is equivalent to the choice of an end-proper ray 
in $M_i$. 

\proclaim{Theorem 5.1} Let $\{M_i\}$ be a countable collection (with 
at least two elements) of connected, oriented, \irr\ open \tms\ each 
of which has only finitely many ends. Suppose $\Gamma$ is a countable,  
locally finite tree whose vertices $v_i$ correspond bijectively to 
the $M_i$. Suppose that to each edge $e_j$ of $\Gamma$ incident with 
$v_i$ we have associated an end-proper ray $\gamma_{i,j}$ in $M_i$. 
Suppose the $\gamma_{i,j}$ are disjoint and that each end of $M_i$ is 
determined by at least one such ray.  Then we have: 
\roster 
\item"{(a)}" The union $\gamma$ of the $\gamma_{i,j}$ is end-proper 
homotopic to an embedded 1-manifold $\gamma\p$ such that the 
corresponding end sum $M$ along $\Gamma$ is a strong end sum.  
\item"{(b)}" $\gamma\p$ can be chosen so that $M$ admits no homeomorphisms 
which reverse orientation or send one summing plane to another or send 
one summand to another. Moreover, 
if the orientation is changed on any one of the summands then the 
resulting end sum $M^*$ is not homeomorphic to $M$. 
\item"{(c)}" There are uncountably many choices of $\gamma\p$ yielding 
pairwise non-homeo\-mor\-phic such $M$. 
\endroster 
\endproclaim 

This will be deduced from Theorem 4.3 and the following result.  

\proclaim{Proposition 5.2} Let $U$ be a connected, orientable, 
\irr, open \3m\ 
with $\mu < \infty$ ends. For each $1 \leq m \leq \mu$ let 
$1 \leq \nu_m < \infty$. Suppose $\be$ is an end-proper 1-manifold 
in $U$ whose components are rays $\be^{m,p}$, where $1\leq m\leq \mu$, 
$1\leq p \leq \nu_m$, and $\be^{m,p}$ determines the $m^{th}$ end 
of $U$. Then $\be$ is end-proper homotopic to a 1-manifold $\al$ 
having the following properties. 
\roster
\item The \3m\  $V$ obtained by removing the interiors of disjoint 
regular neighborhoods $H^{m,p}$ of the $\al^{m,p}$ is \irr, strongly 
aplanar, and \ann\ at $\infty$.  
\item If $\widehat{V}$ is formed by re-attaching, for each end, some, 
but not all, of the halfspaces $H^{m,p}$ to that end, then $\widehat{V}$ 
has all the properties listed in (1). 
\item Each $\widehat{V}$ admits no homeomorphisms which reverse orientation or take 
one component of $\bd \widehat{V}$ to another; distinct $\widehat{V}$ 
are non-homeomorphic. 
\item There are uncountably many choices of the $\al$ yielding 
pairwise non-home\-o\-mor\-phic $V$ with properties (1), (2), and (3); 
moreover the $\widehat{V}$ are pairwise non-homeomor\-phic. 
\endroster
\endproclaim

\demo{Proof of Theorem 5.1} We apply Proposition 5.2 to $M_i$. 
We get strength from statement (1) along with statement (2), 
 which implies that there are no bad summands and hence by Corollary 
3.3 that the sum is non-degenerate. (4) allows us to choose distinct 
$V_i$ to be non-homeomorphic; there are uncountably many such choices 
for the set of all $V_i$. By Theorem 4.3 any homeomorphism between two 
such sums or a sum and itself can be isotoped so as to carry one set 
of summing planes to the other, from which it then follows that it 
must carry summand to corresponding summand. The remainder of the theorem 
then follows from (3) and (4). \qed \enddemo 

The existence of an $\al$ satisfying 
(1)--(4) and having the appropriate distribution of its ends among 
the ends of $U$ was proven in Theorems 6.1, 6.5, and 6.8 of \cite{6}. We 
will briefly outline that construction, modifying it so that $\al$ is 
end-proper homotopic to $\be$. We first list 
some technical tools that we shall use. 

\proclaim{Proposition 5.3} Let $V$ be a connected, \irr, orientable, 
non-compact \3m\ which has a finite number $\mu$ of ends and whose boundary 
consists of a finite number of disjoint planes. Suppose $V$ has an exhaustion 
$\{C_n\}$ such that $C_n \cap \bd V$ consists of a single disk in each 
component of $\bd V$, $C_{n+1}-\Inte C_n$ is \irr, \birr, and \ann, each 
component of $Fr \, C_n$ has negative Euler characteristic and positive 
genus, and $V-\Inte C_n$ has $\mu$ components for all $n \geq 0$. Then 
$V$ is strongly aplanar and \ann\ at $\infty$. \endproclaim 

\demo{Proof} This follows from Lemma 1.3, Theorem 3.4, and Theorem 5.3 of 
\cite{6}. \qed \enddemo 

A compact, connected, \3m $X$ which is not a 3-ball 
is called {\bf excellent} if it is ${\bold P}^2$-\irr\ and \birr,  
contains a 2-sided, proper, \inc\ surface, and every connected, 
proper, \inc\ surface of zero Euler characteristic in $X$ is \bparallel. 
The closure of the complement of a regular neighborhood of a submanifold 
$A$ of a manifold $Q$ is called the {\bf exterior} of $A$ in $Q$. 
A proper 1-manifold $\la$ in a compact \3m $Q$ is called 
{\bf excellent} if its exterior in $Q$ is excellent. We say that 
$\la$ is {\bf poly-excellent} if every non-empty union of the components 
of $\la$ is excellent.

\proclaim{Proposition 5.4} Let $Q$ be a compact, connected, orientable 
\3m\ whose boundary is non-empty and contains no 2-spheres. Suppose 
$\kappa$ is a proper arc in $Q$. The $\kappa$ is homotopic rel $\bd \kappa$ 
to an excellent arc $\lambda$ in $Q$. 
\endproclaim

\demo{Proof} This follows from Theorem 1.1 of \cite{5}. \qed \enddemo 

An {\bf $n$-tangle} $\tau$ is an $n$ component proper 1-manifold 
embedded in a 3-ball such that each component of $\tau$ is an arc. 

\proclaim{Proposition 5.5}  For each $n \geq 1$ poly-excellent 
$n$-tangles exist. \endproclaim

\demo{Proof} This is Theorem 6.3 of \cite{6}. \qed \enddemo 

\proclaim{Lemma 5.6} Let $R$ be a compact, connected \3m. Let $S$ be a 
compact, proper, 2-sided surface in $R$. Let $R\p$ be the \3m obtained 
by splitting $R$ along $S$. Let $S\p$ and $S\pp$ be the two copies of 
$S$ in $\bd R\p$ which are identified to obtain $R$. If each component 
of $R\p$ is excellent, $S\p \cup S\pp$ and $(\bd R\p) - \inte(S\p \cup 
S\pp)$ are \inc\ in $R\p$, and each component of $S$ has negative Euler 
characteristic, the $R$ is excellent. \endproclaim

\demo{Proof} This is Lemma 2.1 of \cite{5}. \qed \enddemo 

\demo{Proof of Proposition 5.2} Let $\{K_n\}$ be an 
exhaustion for $U$. We may assume that each $U-\inte K_n$ 
has $\mu$ components $U^m_n$, $1 \leq m \leq \mu$. Let $Y^m_{n+1}=
U^m_{m+1} \cap (K_{n+1}-\inte K_n)$. By attaching 1-handles to $K_n$ 
inside $U-\inte K_n$ and then passing to a subsequence we may 
assume that $Y^m_{n+1}$ and $G^m_n=K_n \cap Y^m_{n+1}$ are each 
connected and that each $G^m_n$ has genus at least two. Put $\be$ 
in general position 
with respect to $\cup \bd K_n$. We may assume that $\be \cap K_0=
\bd \be$. By attaching 1-handles to $K_n$ whose cores are compact 
components of $\be \cap (U-\inte K_n)$ and passing to a subsequence 
we may further assume that each component of $\be \cap (U-\inte K_n)$ 
is non-compact and thus that $\be^{m,p} \cap (K_{n+1}-\inte K_n)=
\be^{m,p} \cap Y^m_{n+1}$ is an arc $\be^{m,p}_{n+1}$ joining $G^m_n$ 
to $G^m_{n+1}$. Let $D^m_n$ be a disk in $G^m_n$ whose interior 
contains $\be \cap G^m_n$. Let $\be^m_{n+1}=\be^{m,1}_{n+1} \cup 
\cdots \cup \be^{m,\nu_m}_{n+1}$. Let $N^m_{n+1}$ be a regular 
neighborhood of $D^m_n \cup D^m_{n+1} \cup \be^m_{n+1}$ in 
$Y^m_{n+1}$, chosen so that $N^m_{n+1} 
\cap G^m_{n+1}=N^m_{n+2} \cap G^m_{n+2}$. By Proposition 5.4 there 
is an excellent proper arc $\eta^m_{n+1}$ in 
$Y^m_{n+1}-\inte N^m_{n+1}$ such that $\bd \eta^m_{n+1}$ lies in 
$(\bd N^m_{m+1}) \cap \inte Y^m_{n+1}$. Let $T^m_{n+1}$ be the union 
of $N^m_{n+1}$ and a regular neighborhood of $\eta^m_{n+1}$ in 
$Y^m_{n+1}-\inte N^m_{n+1}$. 

$T^m_{n+1}$ is a cube with $\nu_m$ handles whose exterior $L^m_{n+1}$ 
in $Y^m_{m+1}$ is excellent. Moreover, $\be^m_{n+1}$ is a proper 
1-manifold in $T^m_{n+1}$ consisting of unknotted, unlinked arcs. 
We shall homotop $\be^m_{n+1}$ in $T^m_{n+1}$ relative to its 
boundary to obtain a 1-manifold $\al^m_{n+1}$ such that the union 
of the $\al^m_{n+1}$ over all $m$ and $n$ gives the desired union 
of rays $\al$. 

A {\bf classical knot space} $Q$ is a \3m\ \homeo\ to the exterior of a 
non-trivial knot in $S^3$. We say that $Q$ is {\bf incompressibly 
embedded} in a \3m\ $R$ if $Q \subseteq R$ and $\bd Q$ is \inc\ in $R$. 

\proclaim{Lemma 5.7} Let $T$ be a cube with $g$ handles. Let $J_1, \ldots, 
J_{\nu}$ be excellent knots in $S^3$. 
Then there are disjoint classical knot spaces $Q_1, \ldots, 
Q_{\nu}$ in $\inte T$ and disjoint proper arcs $\rho_1, \ldots, 
\rho_{\nu}$ 
in $T-\inte(Q_1 \cup \cdots \cup Q_{\nu})$ such that $Q_p$ is \homeo\ to 
the exterior of $J_p$ in $S^3$, there are disjoint 3-balls $B_p$ in 
$\inte T$ such that $Q_p \subseteq B_p$ and $B_p \cap \rho_q=\emptyset$ 
for $p \neq q$, and for every non-empty subset $\{p_1, \ldots, p_k\}$ of 
$\{1, \ldots, \nu\}$ 
\roster 
\item"{(i)}" the exterior $R$ of the 1-manifold $\rho_{p_1} \cup \cdots \cup 
\rho_{p_k}$ in $T$ is ${\bold P}^2$-\irr, \birr, and \ann, 
\item"{(ii)}" $\bd Q_{p_r}$ is \inc\ in $R$, and 
\item"{(iii)}" given any classical knot space $Q$ incompressibly embedded in 
$\inte R$ there is an ambient isotopy of $Q$ in $R$, fixed on $\bd R$, 
which takes $Q$ to some $Q_{p_r}$. \endroster \endproclaim

\demo{Proof} The case $g=1$ is Lemma 6.7 of \cite{6}. Choose disjoint proper 
disks $D_1, \ldots, D_g$ in $T$ which split it to a 3-ball $B$. 
Let $f:B \rightarrow T$ be the identification map. 
Let $Z_1, \ldots, Z_{\nu}$ be disjoint disks in $\inte 
D_g$. Let $T_p$, $1 \leq p \leq \nu$, be disjoint regular neighborhoods 
of $\bd Z_p$ in $T$, chosen so that $A_p=D_g \cap T_p$ is a regular 
neighborhood of $\bd Z_p$ in $D_g$. Then $f^{-1}(T_p)$ is the union of 
two solid tori $T\p_p$ and $T\pp_p$. Let $B^*=B-Int(\cup_{p=1}^{\nu} 
T\p_p \cup T\pp_p)$. Then $f(\bd B^*) \cap (D_1 \cup \cdots \cup D_g)$ 
consists of $D_1, \ldots, D_{g-1}$, 
together with disks $E_1, \ldots, E_{\nu}$, and a disk with $\nu$ holes 
$P$ contained in $D_g$. Whenever $S$ is one of these surfaces, 
$f^{-1}(S)$ consists of two surfaces $S\p$ and $S\pp$. 
We let $T^*$ denote $T-\inte(T_1 \cup \cdots \cup T_{\nu})$. 

By Proposition 5.5 $B^*$ contains a poly-excellent $2(g+1)\nu$-tangle. 
We denote its components by $\rho_{p,j}$, $1 \leq p \leq \nu$, 
$1 \leq j \leq 2(g+1)$. Isotop this tangle so that $\rho_{p,1}$ runs 
from $\inte f^{-1}(\bd T)$ to $\inte D\p_1$, 
$\rho_{p,2i}$ runs from $\inte D\pp_i$ to itself for $1 \leq i \leq g-1$, 
$\rho_{p,2i+1}$ runs from $\inte D\p_i$ to $\inte D\p_{i+1}$ for 
$1 \leq i \leq g-2$, $\rho_{p,2g-1}$ runs from $\inte D\p_{g-1}$ to 
$\inte E\p_p$, $\rho_{p,2g}$ runs from $\inte E\pp_p$ to itself, 
$\rho_{p,2g+1}$ runs from $\inte E\p_p$ to $P\p$, and $\rho_{p,2g+2}$ 
runs from $P\pp$ to $\inte f^{-1}(\bd T)$. We further require that the 
endpoints of the $\rho_{p,j}$ match up in such a way that 
$f(\cup_{j=1}^{2g+2}\rho_{p,j})$ is a proper arc $\rho_p$ in $T^*$. 

We now glue the exteriors $Q_p$ of the knots 
$J_p$ to $T^*$ so that a meridian of $J_p$ is identified with $\bd E_p$ 
and let $B_p$ be the union of $Q_p$ and a regular neighborhood of $E_p$ 
in $T^*$. Then $B_p$ is a 3-ball, from which it follows that this space 
is again a cube with $g$ handles, which we denote again by $T$. 

Now suppose that we have a non-empty subset $\{p_1, \ldots, p_k\}$ of 
$\{1, \ldots, \nu\}$. It follows from Lemma 5.6 that 
$\rho_{p_1} \cup \cdots \cup \rho_{p_k}$ is excellent in 
$T-\inte (Q_{p_1} \cup \cdots \cup Q_{p_k})$. Standard general position 
and isotopy arguments now show that the exterior $R$ of this 1-manifold 
is ${\bold P}^2$-\irr, \birr, and \ann, that each $\bd Q_{p_r}$ is \inc\ in $R$ 
and that every incompressible torus in $R$ is isotopic to one of these 
tori. The result then follows. \qed \enddemo 

We now complete the proof of Proposition 5.2. Let $\al^m_{n+1}$ be 
a proper 1-manifold in $T^m_{n+1}$ consisting of $\nu_m$ arcs having 
the properties stated in Lemma 5.7. Denote the classical knot spaces 
involved by $Q^{m,p}_{n+1}$. Since $\pi_1(\bd T^m_{n+1}) 
\rightarrow \pi_1(T^m_{n+1})$ is onto we may isotop $\al^m_{n+1}$ so 
that $\bd \al^{m,p}_{n+1}=\bd \be^{m,p}_{n+1}$ and $\al^m_{n+1}$ and 
$\be^m_{n+1}$ are homotopic relative to this common boundary. Thus 
the union $\al$ of the $\al^{m,p}_{n+1}$ is end-proper homotopic 
to $\be$. 

Now for each $m$, $1 \leq m \leq \mu$, choose a non-empty subset 
of $\{1, \ldots, \nu_m\}$. 
By property (i) of Lemma 5.7 the exterior $R^m_{n+1}$ in $T^m_{n+1}$ 
of the corresponding  
union of components of $\al^m_{n+1}$ is \irr, \birr, and \ann. Since 
the same is true of the exterior $L^m_{n+1}$ of $T^m_{n+1}$ in 
$Y^m_{n+1}$ standard general position and isotopy arguments show that 
these properties hold for $X^m_{n+1}=L^m_{n+1} \cup R^m_{n+1}$. 
Let $X_{n+1}=X^1_{n+1} \cup \cdots \cup X^{\mu}_{n+1}$, 
$C_0=K_0$, and $C_{n+1}=K_0 \cup X_1 \cup \cdots \cup X_{n+1}$. 
Then $\{C_n\}$ is an exhaustion for the exterior $\widehat{V}$ of 
the corresponding components of $\al$ in $U$. The application of 
Proposition 5.3 to $\{C_n\}$ now implies properties (1) and (2) of 
Proposition 5.2. 

The proof of properties (3) and (4) of Proposition 5.2 is identical 
to that of Theorem 6.8 of \cite{6}, with Lemma 6.7 of that paper replaced 
by Lemma 5.7 of this paper. For the sake of completeness we briefly 
recall the construction, referring the reader to \cite{6} for details. 
We choose a countably infinite family of excellent knots in $S^3$ 
whose exteriors admit no orientation reversing homeomorphisms. 
We index this family by quadruples $(m,p,n,q)$, where $1\leq m\leq \mu$, 
$1\leq p \leq \nu_m$, $n\geq 0$, and $q \in \{0,1\}$. We choose some 
function $q=\varphi(m,p,n)$ and then carry out our construction with 
the knot space $Q(m,p,n,q)$ associated to the arc $\al^{m,p}_{n+1}$ 
as in Lemma 5.7; this produces a manifold $V[\varphi]$. Given a  
collection $E^{m,p}$ of boundary planes which includes at least one 
plane from each end and any compact subset of 
$V[\varphi]$ meeting exactly these boundary planes, it turns out that there is 
a larger compact subset meeting exactly these boundary planes such that the 
incompressibly embedded knot spaces in the complement of this subset 
are precisely the corresponding $Q(m,p,n,q)$ which it contains. This 
property, together with the fact that the knot spaces admit no 
orientation reversing homeomorphisms and the fact that the set of 
functions $\varphi$ is uncountable, yields properties (3) and (4). 
\qed \enddemo

\proclaim{Theorem 5.8} There are uncountably many connected, \irr, open 
3-man\-i\-folds which are end-prime but not \rirr. \endproclaim 

\demo{Proof} Let $U$ be any connected, oriented, \irr, open \3m\ 
with one end. By Proposition 5.2 we may choose two disjoint proper 
halfspaces in $U$ such that the \3m\ $V$ obtained by removing their 
interiors is \irr, strongly aplanar, and anannular at infinity. Let 
$M$ be the \3m\ obtained by gluing the two components of $\bd V$ together 
via an orientation reversing homeomorphism. Since the plane $E$ in $M$ 
which is the image of $\bd V$ is non-separating $M$ cannot be \rirr. 
Suppose it were not end-prime. Then $M$ would be the plane sum of 
\tms\ $V_1$ and $V_2$ each having boundary a plane but neither being 
\homeo\ to $\halfspace$. Thus the plane $P$ along which the sum is 
taken is non-trivial in $M$. Then by Theorem 4.1 we have that $P$ is 
ambient isotopic to $E$, which cannot happen since $P$ separates $M$. 
\qed \enddemo 

\head 6. End-prime Decompositions \endhead

In this section we show that there are \irr\ open  \tms\ which are not 
end-prime and do not admit decompositions via end sum into end-prime 
\tms. We first show that one can reduce to the case of non-degenerate 
end sum decompostions. 

\proclaim{Theorem 6.1} Suppose $M$ is an end sum of connected, \irr, 
end-prime open \tms\ $M_i$  along a tree $\Gamma$. Then either $M$ is 
end-prime or can be expressed as a non-degenerate end sum of connected, 
\irr, end-prime open \tms. \endproclaim

\demo{Proof} Let $\{V_i\}$ and $\{E_j\}$ denote the vertices and edges of 
the corresponding plane sum. 

Suppose some of the summing planes are trivial. By Theorem 3.2 each such 
plane determines a bad branch. Partially order the bad branches by 
inclusion. If there are no maximal bad branches, then $M$ is the monotone 
union of a sequence $\{\be_j\}$ of bad branches associated with a sequence 
$\{E_j\}$ of summing planes. Since the union of all summing planes is proper 
in $M$ it follows that every compact subset of $M$ lies in the interior 
of some branch, and hence $M$ is homeomorphic to $\RRR$ and so is end-prime. 

Thus we may assume that maximal bad branches exist. Suppose two such 
branches $\be_j$ and $\be_k$ intersect. If $\be_j \cap \be_k$ is a 
summing plane, then again $M$ is \homeo\ to $\RRR$. So we may assume that 
$E_j$ lies in the interior of $\be_k$. Since $\be_k$ is \homeo\ to 
$\halfspace$, which is \rirr, $E_j$ is trivial in $\be_k$, from which it 
follows that $M$ is \homeo\ to $\RRR$. 

Hence we may assume that distinct maximal bad branches are disjoint. 
All the trivial summing planes are contained in the union of the maximal 
bad branches. Suppose $V_i$ is a plane summand which is not contained in 
a bad branch. If $V_i$ meets a bad branch $\be_j$ in a summing plane $E_j$, 
then $\be_j$ is maximal. Let $V_i\p$ be the union of $V_i$ with all such 
$\be_j$. Then $V_i$ and $V_i\p$ each have interior \homeo\ to $M_i$. 
We can thus eliminate all those $M_k$ whose corresponding $V_k$ lie in a 
bad branch. If only one summand remains, then $M$ is end-prime. If more 
than one summand remains, the $M$ is expressed as an end sum of end-prime 
\tms\ in which no summing plane is trivial. 

Now suppose some pairs of distinct summing planes are parallel. 
By Theorem 3.2 they determine a bad trail. Since there are no trivial 
summing planes there are no bad branches and hence this trail has no side 
branches. Each of its vertices is therefore \homeo\ to $\prodspace$. 
Partially order the bad trails by inclusion. If there are no maximal elements 
then there is an infinite nested sequence of bad trails in $M$ and 
therefore a trivial summing plane. Thus maximal bad trails exist and 
clearly distinct such arcs are disjoint. Fix a summand $V_0$ which is 
not contained in a bad trail. Let $V_0\p$ be the union of $V_0$ with any 
bad trails which meet it. Both $V_0$ and $V_0\p$ have interior \homeo\ 
to $M_0$, so we can delete all the $M_k$ whose corresponding $V_k$ are 
contained in these bad trails. We now repeat this argument with each of the 
$V_i$ which meet $V_0\p$ in a summing plane and then continue to 
inductively define a new end sum in this manner. The set of end 
summands is a subset of the original set and is clearly non-degenerate. 
\qed \enddemo 

A connected, non-compact \3m\ is {\bf eventually end-\irr\ } if it admits an 
exhaustion $\{C_n\}$ such that $Fr \, C_n \cup Fr \, C_{n+1}$ is 
\inc\ in $C_{n+1}-\Inte C_n$ for all $n \geq 0$. 

\proclaim{Theorem 6.2} There exists an \irr, eventually end-\irr, 
contractible open \3m\ $M$ which is not end-prime and does not admit 
a decomposition via end sum into end-prime \tms. \endproclaim 

\demo{Proof} By Theorems 4.3 and 6.1 it suffices to construct an \irr, 
eventually end-\irr, contractible open \3m\ $M$ which is \homeo\ to a 
strong end sum of itself with itself. We shall construct such an $M$ 
having an exhaustion $\{C_n\}$ where $C_n$ is a cube with $2^{n+1}$ 
handles. We shall present $M$ as the direct limit of a sequence 
$K_0 @> g_0 >> K_1 @> g_1 >> K_2 @> g_2 >> \cdots$, 
where $K_n$ is a cube with 
$2^{n+1}$ handles and $g_n$ embeds $K_n$ into the interior of $K_{n+1}$. 
The image of $K_n$ in the direct limit will be $C_n$. 

We regard $K_n$ as two copies $W^{\pm}_n$ of a cube with $2^n$ handles 
joined by a 1-handle $L_n=D_n \times [-1,1]$, where $D_n \times \{\pm 1\}$ 
are disks lying in $\bd W^{\pm}_n$, respectively. We denote $D_n \times 
\{0\}$ by $D_n$ and the center of $D_n$ by $x_n$. We let $L^-_n$ and 
$L^+_n$ denote $D_n \times [-1,0]$ and $D_n \times [0,1]$, $D^{\pm}_n$ 
denote $D_n \times \{\pm 1\}$, and $x^{\pm}_n$ denote $\{x_n\} \times 
\{\pm 1\}$, respectively. We let $p_n:D_n \times [-1,1] \rightarrow D_n$ 
be projection onto the first factor. Let $r_n$ be an orientation preserving 
involution of $K_n$ which interchanges $W^+_n$ and $W^-_n$ as well as 
$L^+_n$ and $L^-_n$; we assume that the restriction of $r_n$ to $L_n$ has 
the form $r_n(x,t)=(s_n(x),-t)$, where $s_n$ is reflection in a diameter 
of $D_n$. Thus $r_n$ fixes $x_n$ and interchanges $x^+_n$ and $x^-_n$. 
We also choose an orientation preserving homeomorphism $f_n:K_n 
\rightarrow W^+_{n+1}$. We next describe the embeddings $g_n:K_n 
\rightarrow K_{n+1}$. 

Let $h_0:W^+_0 \rightarrow \inte W^+_1$ be any null-homotopic embedding. 
We define $g_0$ to be $h_0$ on $W^+_0$. By Proposition 5.4 there is 
an excellent arc $\al_0$ in $W^+_1-\inte g_0(W^+_0)$ joining $x^+_1$ to 
$\bd g_0(W^+_0)$. We extend $g_0$ over $D_0 \times [\frac{1}{2},1]$ by 
sending it to a regular neighborhood $N_0$ of $\al_0$ in 
$W^+_1-\inte g_0(W^+_0)$. We require that $g_0$ take $D_0 \times \{\frac
{1}{2}\}$ to a disk in the interior of $D^+_1$ whose image under $p_1$ is 
invariant under $r_1$. We then extend $g_0$ over $D_0 \times 
[0,\frac{1}{2}]$ by setting $g_0(x,t)=(p_1(g_0(x,\frac{1}{2})),2t)$. 
Finally, we extend $g_0$ over $W^-_0 \cup L^-_0$ by defining it to be 
$r_1g_0r_0$. Thus we have $g_0:K_0 \rightarrow K_1$. Note that the image 
is invariant under $r_1$ and that $r_1g_0=g_0r_0$. 

We now define $g_1:K_1\rightarrow K_2$. The basic idea is to embed 
$W^+_1$ in $W^+_2$ in the same fashion in which $K_0$ is embedded in $K_1$, 
then choose an excellent arc joining $x^+_2$ to the boundary of the image 
of this embedding, and then extend the embedding to all of $K_1$ as in 
the previous step. More precisely, we let $h_1=f_1g_0f^{-1}_0: W^+_1 
\rightarrow W^+_2$ be our initial embedding. We define $g_1$ on $W^+_1$ 
to be $h_1$. Let $\al_1$ be an excellent arc in $W^+_2-\inte g_1(W^+_1)$ joining 
$x^+_2$ to $\bd g_1(W^+_1)$. We extend $g_1$ over $D_1 \times 
[\frac{1}{2}, 1]$ by sending it to a regular neighborhood $N_1$ of $\al_1$ 
in $W^+_2-\inte g_1(W^+_1)$, requiring that $p_2(g_1(D_1 \times 
\{\frac{1}{2}\}))$ be invariant under $r_2$. We extend $g_1$ over 
$D_1 \times [0,\frac{1}{2}]$ by setting $g_1(x,t)=(p_2(g_1(x,\frac{1}{2})), 
2t)$. Finally, we extend $g_1$ over $W^-_1 \cup L^-_1$ by defining it to 
be $r_2g_1r_1$. Thus we have $g_1:K_1 \rightarrow K_2$, an embedding whose 
image is invariant under $r_2$ such that $r_2g_1=g_1r_1$. 

Continuing in this manner we construct sequences of embeddings 
$g_n:K_n\rightarrow K_{n+1}$, $h_n:W^+_n\rightarrow W^+_{n+1}$, 
where $h_n=f_ng_{n-1}f_{n-1}$, $g_n\_{W^+_n}=h_n$, $g_n(D_n \times 
[\frac{1}{2},1])$ is a regular neighborhood $N_n$ of an excellent arc 
$\al_n$ joining $x^+_n$ to $\bd g_n(W^+_n)$, $g_n(x,t)=
(p_{n+1}(g_n(x,\frac{1}{2})),2t)$ for $(x,t) \in D_n \times 
[0,\frac{1}{2}]$, and letting $g_n$ be $r_{n+1}g_nr_n$ on 
$W^-_n \cup L^-_n$. 
We then let $M$ be the direct limit of the sequence 
$K_0 @> g_0 >> K_1 @> g_1 >> K_2 @> g_2 >> \cdots$ and denote the image of 
$K_n$ in $M$ by $C_n$. We let $V_n$ be the image of $W^+_n$ in $M$, 
and let $V$ be the union of the $V_n$. Note that the $r_n$ induce 
an involution $r$ of $M$, $M=V \cup r(V)$, and $V \cap r(V)=\bd V= 
\bd r(V)$, which is a plane $E$ invariant under $r$. 

\proclaim{Lemma 6.3} $M$ is \irr, eventually end-\irr, and contractible. 
$V$ is \irr, strongly aplanar, and \ann\ at $\infty$. \endproclaim 

\demo{Proof} $M$ is \irr\ because it is a monotone union of cubes with 
handles. $C_{n+1}-\inte C_n$ is \homeo\ to $K_{n+1}-\inte g_n(K_n)$, 
which is the union of $W^+_{n+1}-\Inte (g_n(W^+_n) \cup N_n)$ and its 
image under $r_{n+1}$. Since $\al_n$ is an excellent arc these manifolds are 
\irr, \birr, and \ann. They intersect in an \inc\ annulus, and so by 
a standard general position and isotopy argument their union is \irr\ 
and \birr. Since $h_0:W^+_0\rightarrow 
W^+_1$ is null-homotopic we have that $g_0:K_0\rightarrow K_1$ is 
null-homotopic, and so $h_1=f_1g_0f^{-1}_0:W^+_1\rightarrow W^+_2$ is 
null-homotopic. Thus by induction we have that all the $h_n$ and $g_n$ 
are null-homotopic, and so $M$ is contractible. 

The irreducibility of $V$ follows from that of $M$ and the fact that 
$E$ is a plane. The strong aplanarity and anannularity at $\infty$ of 
$V$ follow from Proposition 5.3 and the anannularity of 
$W^+_{n+1}-\Inte (g_n(W^+_n)\cup N_n)$. \qed \enddemo 

\proclaim{Lemma 6.4} $\Inte V$ is \homeo\ to $M$. \endproclaim 

\demo{Proof} $\Inte V$ is the direct limit of the sequence 
$W^+_1 @> h_1 >> W^+_2 @> h_2 >>  
W^+_3 @> h_3 >> \cdots$. Since the homeomorphisms 
$f_n:K_n\rightarrow W^+_{n+1}$ satisfy $f_{n+1}g_n=h_nf_n$ we have 
that they induce a homeomorphism $f:M \rightarrow \Inte V$. \qed \enddemo 

This completes the proof of Theorem 6.2. \qed \enddemo

\enddocument